\documentclass[12pt]{article}
\usepackage[none]{hyphenat}

\usepackage{amssymb,amsmath,amsfonts,amsthm,amsxtra}
\usepackage{fancyhdr, lastpage}

\usepackage[utf8]{inputenc}

\newcommand\largediamond{\mbox{\large$\diamondsuit$}}
\newcommand{\mr}{\mathrm}
\newcommand{\mc}{\mathcal}

\newcommand{\op}{\operatorname}
\newcommand{\boldsigma}{\mathbf{\mathop{\pmb{\sum}}}}

\newtheorem{theorem}{Theorem}[section]

\newenvironment{ackn}{\begin{trivlist} \item[] {\bf Acknowledgement.}}{\hspace*{0pt}\end{trivlist}}
\newenvironment{con}{\begin{trivlist} \item[] {\bf Conjecture.}}{\hspace*{0pt}\end{trivlist}}
\newenvironment{Proof}{\begin{trivlist} \item[] {\bf Proof.}}{\hspace*{0pt}\hfill\qedsymbol\end{trivlist}}

\newenvironment{Proofs}{\begin{trivlist} \item[] {\bf Proofs.}}{\hspace*{0pt}\hfill\end{trivlist}}
\newenvironment{PROOF}{\begin{trivlist} \item[] {\bf Proof.}}{\end{trivlist}}

\newsavebox{\Prfref}
\newenvironment{Prf}[1][{\bf.}]
{\sbox{\Prfref}{\textbf{#1.}}
\begin{trivlist} \item[] {\bf Proof of Theorem \usebox{\Prfref}}}{\hspace*{0pt}\hfill $\Box$\end{trivlist}}

\newsavebox{\prfref}
\newenvironment{prf}[1][{\bf.}]
{\sbox{\prfref}{\textbf{#1.}}
\begin{trivlist} \item[] {\bf Proof of Theorem \usebox{\prfref}}}{\end{trivlist}}

\newsavebox{\prfcorref}
\newenvironment{proofcor}[1][{\bf.}]
{\sbox{\prfcorref}{\textbf{#1.}}
\begin{trivlist} \item[] {\bf Proof of Corollary \usebox{\prfcorref}}}{\hspace*{0pt}\hfill $\Box$\end{trivlist}}

\newtheoremstyle{break}% an optional reference is inserted after the theorem number
{\topsep}%	Space above
{\topsep}%	Space below
{\it}%         		Body font
{}%			Indent amount (empty = no indent, \parindent = para indent)
{}%			Thm head font
{}%        		Punctuation after thm head
{ }%	Space after thm head: " " = normal interword space;	\newline = linebreak
{\thmname{\textbf{#1}}\thmnumber{ \textbf{#2.}}}% Thm head spec (can be left empty, meaning `normal')

\theoremstyle{break}

\newtheoremstyle{ref}% an optional reference is inserted after the theorem number
{\topsep}	%      Space above
{\topsep}	%      Space below
{\it}%         		Body font
{}%         		Indent amount (empty = no indent, \parindent = para indent)
{}%			Thm head font
{}%        		Punctuation after thm head
{ }%			Space after thm head: " " = normal interword space;
{\thmname{{\bfseries#1}}\thmnumber{ \textbf{#2\thmnote{\rm #3}\textbf .}}}%Thm head spec (can be left empty, meaning `normal')

\theoremstyle{ref}
\newtheorem{lem}[theorem]{Lemma}
\newtheorem{prop}[theorem]{Proposition}
\newtheorem{cor}[theorem]{Corollary}
\newtheorem{thm}[theorem]{Theorem}

\newtheorem*{thmm}{Main Theorem}

\newtheoremstyle{nnref}% an optional reference is inserted after the theorem, there is no theorem number
{\topsep}	%      Space above
{\topsep}	%      Space below
{}%         		Body font
{}%         		Indent amount (empty = no indent, \parindent = para indent)
{}%			Thm head font
{}%        		Punctuation after thm head
{ }%			Space after thm head: " " = normal interword space;
{\thmname{\textbf{#1}\thmnote{\textrm{ #3}}\textbf{.}}}%Thm head spec (can be left empty, meaning `normal')

\theoremstyle{nnref}
\newtheorem{defn}{Definition}
\newtheorem{rem}{Remark}

\newtheoremstyle{prime}% an apostrophe is inserted after the theorem number
{\topsep}	%      Space above
{\topsep}	%      Space below
{\it}%         		Body font
{}%         		Indent amount (empty = no indent, \parindent = para indent)
{}%			Thm head font
{}%        		Punctuation after thm head
{ }%			Space after thm head: " " = normal interword space;
{\thmname{{\bfseries#1}}\thmnumber{ \textbf{#2$\boldsymbol{'}$\thmnote{\rm #3}\textbf .}}}%Thm head spec (can be left empty, meaning `normal')

\theoremstyle{prime}
\newtheorem{thmp}[theorem]{Theorem}
\newtheorem{corp}[theorem]{Corollary}

\begin{document}
%\pagestyle{fancyplain}
%\fancyhf{}
%\rhead{Franklin D. Tall 4872}
%\lhead{Page \thepage\space of \pageref{LastPage}}
\sloppy
%\begin{titlepage}
\title{PFA$(S)[S]$ and Locally Compact Normal Spaces}
\author{Franklin D. Tall\makebox[0cm][l]{$^1$}}

\footnotetext[1]{Research supported by NSERC grant A-7354.\vspace*{2pt}}
\date{\today}
\maketitle

\begin{abstract}
We examine locally compact normal spaces in models of form PFA$(S)[S]$, in particular characterizing paracompact, countably tight ones as those which include no perfect pre-image of $\omega_1$ and in which all separable closed subspaces are Lindel\"{o}f.
\end{abstract}

\renewcommand{\thefootnote}{}
\footnote
{\parbox[1.8em]{\linewidth}{$2010$ AMS Math.\ Subj.\ Class.\ Primary 54A35, 54D15, 54D20, 54D45, 03E35, 03E65; Secondary 54E35.}\vspace*{5pt}}
\renewcommand{\thefootnote}{}
\footnote
{\parbox[1.8em]{\linewidth}{Key words and phrases: $\mr{PFA}(S)[S]$, paracompact, locally compact, normal, perfect pre-image of $\omega_1$, locally connected, reflection, Axiom R, P-ideal Dichotomy, Dowker space, collectionwise Hausdorff, homogeneous compacta.}}

\section{Introduction}
We will be using a particular kind of model of set theory constructed with the aid of a supercompact cardinal. These models have been used in~\cite{Larson,Larsona,TallTopCons,TallDowker,T2,T,Todorcevic}. We start with a particular kind of Souslin tree --- a \textit{coherent} one --- in the ground model. Such trees are obtainable from $\largediamond$~\cite{Larson1999,Shelah1999}. Their definition will not concern us here. One then iterates proper partial orders as in the proof of the consistency of $\mr{PFA}$ (so we need to assume the consistency of a supercompact cardinal) but only those that preserve the Souslinity of that tree. By~\cite{Miyamoto1993}, that produces a model for \textit{$\mathit{PFA}(S)$: $\mr{PFA}$ restricted to partial orders that preserve $S$}. We then force with $S$. We shall say \textit{$\mathit{PFA}(S)[S]$ implies $\varphi$} if $\varphi$ holds whenever we force with $S$ over a model of $\mr{PFA}(S)$, for $S$ a coherent Souslin tree. We shall say \textit{$\varphi$ holds in a model of form $\mathit{PFA}(S)[S]$} if $\varphi$ holds in some particular model obtained this way.

PFA$(S)[S]$ and particular models of it impose a great deal of structure on locally compact normal spaces because they entail many useful consequences of both PFA and $V = L$.  We amalgamate here three previous preprints \cite{T1}, \cite{TallDowker}, and \cite{T2} dealing with characterizing paracompactness and killing Dowker spaces in locally compact normal spaces, as well as with homogeneity in compact hereditarily normal spaces.  Our proofs will avoid the difficult set-theoretic arguments in other papers on PFA$(S)[S]$ by just quoting the familiar principles derived there, and so should be accessible to any set-theoretic topologist.

The consequences of PFA$(S)[S]$ we shall use and the references in which they are proved are:

\begin{description}
\item [(Balogh's) $\boldsigma$] (defined below) \cite{fischerTallNew};
\item [$\aleph_1$-CWH] (locally compact normal spaces are $\aleph_1$-collectionwise Hausdorff \cite{TallTopCons});
\item [PID] (P-ideal Dichotomy (defined below) \cite{LarsonErice}).
\end{description}
We also mention for the reader's interest:
\begin{description}
\item [MM] (compact countably tight spaces are sequential \cite{Todorcevic2008});
\item [FCL] (every first countable hereditarily Lindel\"{o}f space is hereditarily separable \cite{LT});
\item [FC$\aleph_1$-CWH] (every first countable normal space is $\aleph_1$-collectionwise Hausdorff \cite{Larson});
\item [OCA] (Open Colouring Axiom \cite{F});
\item [$\mathfrak{b} = \aleph_2$] (\cite{Larson1999}).
\end{description}
In the particular model used in \cite{Larson}, we also have:
\begin{description}
\item [FCCWH] (every first countable normal space is collectionwise Hausdorff \cite{Larson});
\item [CWH] (locally compact normal spaces are collectionwise Hausdorff \cite{TallTopCons});
\item [Axiom R] (see below \cite{Larsona}).
\end{description}
Proofs of all of these results are published or available in preprints, with the exception of \textbf{MM}, for which only a brief sketch exists.  We will therefore avoid using it.

In Section 2 we characterize paracompactness in locally compact normal spaces in certain models of PFA$(S)[S]$.  In Section 3, we improve our characterization via the use of P-ideal Dichotomy.  In Section 4, we examine applications of Axiom R.  In Section 5 we obtain some reflection results in ZFC for (locally) connected spaces.  In Section 6, we apply our results to exclude certain locally compact Dowker spaces in models of PFA$(S)[S]$.  In Section 7, we apply PFA$(S)[S]$ to compact hereditarily normal spaces.

\section{Characterizing paracompactness in locally compact normal spaces}
Engelking and Lutzer~\cite{Engelking1977} characterized paracompactness in generalized ordered spaces by the absence of closed subspaces homeomorphic to stationary subsets of regular cardinals. This was extended to monotonically normal spaces by Balogh and Rudin~\cite{Balogh1992}. Moreover, for first countable generalized ordered spaces, one can do better:

\begin{prop}[~\cite{Tall1994}]
Assuming the consistency of two supercompact cardinals, it is consistent that a first countable generalized ordered space is (hereditarily) paracompact if and only if no closed subspace of it is homeomorphic to a stationary subset of $\omega_1$.
\end{prop}

We were interested in consistently obtaining a similar characterization for locally compact normal spaces. However, as we shall see, the locally compact, separable, normal, first countable, submetrizable, non-paracompact space Weiss constructed in~\cite{Weiss1980} has no subspace homeomorphic to a stationary subset of $\omega_1$. Nonetheless, for restricted classes of locally compact normal spaces, we can get characterizations of paracompactness that do depend on the spaces' relationship with $\aleph_1$.

We will assume all spaces are Hausdorff.

In~\cite{Larsona} we proved:
\begin{lem}[]
\label{lchncounter}
In a particular model (the one of \cite{Larson}) of form $\mr{PFA}(S)[S]$, every locally compact, hereditarily normal space which does not include a perfect pre-image of $\omega_1$ is paracompact.
\end{lem}

Lemma~\ref{lchncounter} will follow from what we prove here as well. We can turn this result into a characterization as follows.

\begin{thm}
\label{lchnshp}
There is a model of form $\mr{PFA}(S)[S]$ in which locally compact hereditarily normal spaces are (hereditarily) paracompact if and only if they do not include a perfect pre-image of $\omega_1$.
\end{thm}
\begin{Proof}
The backward direction follows from Lemma~\ref{lchncounter}, since a space is hereditarily paracompact if every open subspace of it is paracompact, and open subspaces of locally compact spaces are locally compact. The ``hereditarily'' version of the other direction is because perfect pre-images of $\omega_1$ are countably compact and not compact, and hence not paracompact.  Without ``hereditarily'' we need:

\begin{lem}[~\cite{EisworthNyikos}]\label{lem2p4}
In a countably tight space, perfect pre-images of $\omega_1$ are closed.
\end{lem}

\begin{lem}[~\cite{Balogh1983,Balogh2002,Larsona}]
A locally compact space has a countably tight one-point compactification if and only if it does not include a perfect pre-image of $\omega_1$.
\end{lem}

Note that if $X$ has a countably tight one-point compactification, $X$ itself is countably tight.
\end{Proof}

We also now have a partial characterization for locally compact spaces that are only normal:

\begin{thm}\label{thm2p6}
There is a model of form PFA$(S)[S]$ in which a locally compact normal space is paracompact and countably tight if and only if its separable closed subspaces are Lindel\"{o}f and it does not include a perfect pre-image of $\omega_1$.
\end{thm}

The proof of Theorem \ref{thm2p6} is quite long.  It is convenient to first prove the weaker

\begin{thm}
\label{mfpfa}
There is a model of form $\mr{PFA}(S)[S]$ in which a locally compact normal space $X$ is paracompact and countably tight if and only if the closure of every Lindel\"{o}f subspace of $X$ is Lindel\"{o}f, and $X$ does not include a perfect pre-image of $\omega_1$.
\end{thm}

One direction is easy. As we saw earlier, perfect pre-images of
$\omega_1$ will be excluded by countable tightness plus
paracompactness. It is also easy to see that a paracompact space with
a dense Lindel\"{o}f subspace is Lindel\"{o}f --- since it has
countable extent --- so closures of Lindel\"{o}f subspaces are
Lindel\"{o}f. The other direction is harder, but much of the work has
been done elsewhere. We refer to~\cite{Fleissner1986} for a definition
of the reflection axiom \textit{Axiom $R$}. Dow~\cite{Dow1992} proved
it equivalent to stationary set reflection for stationary subsets of
$[\kappa]^{\omega}$, $\kappa$ uncountable. However, we shall only use
the following three results concerning it. We have:
\begin{lem}[~\cite{Larsona}]
Axiom R holds in the PFA$(S)[S]$ model of \cite{Larson}.
\end{lem}

\begin{defn}
  $L(Y)$, the Lindelöf number of $Y$, is the least cardinal $\kappa$ such
  that every open cover of $Y$ has a subcover of size $\le \kappa$.
\end{defn}

\begin{lem}[~\cite{Balogh2002}]
Axiom $R$ implies that if $X$ is a locally Lindel\"{o}f, regular, countably tight space such that every open $Y$ with $L(Y)\leq\aleph_1$ has $L(\overline{Y})\leq\aleph_1$, then if $X$ is not paracompact, it has a clopen non-paracompact subspace $Z$ with $L(Z)\leq\aleph_1$.
\end{lem}

\begin{lem}[~\cite{Balogh2002}]
\label{rimp}
Axiom $R$ implies that if $X$ is locally Lindel\"{o}f, regular, countably tight, and not paracompact, then $X$ has an open subspace $Y$ with $L(Y)\leq\aleph_1$, such that $Y$ is not paracompact.
\end{lem}

We also have:
\begin{lem}[]
\label{yslls}
If $Y$ is a subset of a locally Lindel\"{o}f space of countable tightness in which closures of Lindel\"{o}f subspaces are Lindel\"{o}f, then if $L(Y)\leq\aleph_1$, then $L(\overline{Y})\leq\aleph_1$.
\end{lem}
\begin{Proof}
For let $\mc{U}$ be a collection of open sets with Lindel\"{o}f closures covering $\overline{Y}$. There are $\aleph_1$ of them, say $\{ U_{\alpha}\}_{\alpha<\omega_1}$, which cover $Y$. Then $\overline{\bigcup_{\alpha<\omega_1}U_{\alpha}} = \bigcup_{\alpha<\omega_1}\overline{\bigcup_{\beta<\alpha}U_{\beta}} = \bigcup_{\alpha < \omega_1}\overline{\bigcup_{\beta < \alpha}\overline{U}_{\beta}} \supseteq\overline{Y}$. But each $\overline{\bigcup_{\beta<\alpha}\overline{U}_{\beta}}$ is Lindel\"{o}f, so $L\big(\overline{\bigcup_{\alpha<\omega_1}U_{\alpha}}\big)\leq\aleph_1$. $\overline{Y}$ is a closed subspace of $\overline{\bigcup_{\alpha<\omega_1}U_{\alpha}}$, so it too has Lindel\"{o}f number $\leq\aleph_1$.
\end{Proof}

To finish the proof of Theorem~\ref{mfpfa} it therefore suffices to prove:
\begin{thm}
\label{nppthm}
$\mr{PFA}(S)[S]$ implies that if $X$ is a locally compact normal space
with $L(X)\leq\aleph_1$, closures of Lindel\"{o}f subspaces of $X$ are
Lindel\"{o}f, and $X$ includes no perfect pre-image of $\omega_1$,
then $X$ is paracompact.
\end{thm}

Crucial ingredients in proving this are \textbf{$\aleph_1$-CWH} and:

\begin{lem}[~\cite{fischerTallNew}, \cite{Todorcevic}]
\label{todcount}
$\mr{PFA}(S)[S]$ implies
\begin{description}
\item [$\boldsigma$:] {if $X$ is compact and countably tight, and $Z \subseteq X$ is such that $|Z| \leq \aleph_1$ and there exists a collection $\mc{V}$ of open sets, $|\mc{V}| \leq \aleph_1$, and a collection $\mc{U} = \{U_V: V \in \mc{V}\}$ of open sets, such that $Z \subseteq \bigcup \mc{V}$, and for each $V \in \mc{V}$, there is a $U_V \in \mc{U}$ such that $V \subseteq \overline{V} \subseteq U_V$, and $|U_V \cap Z| \leq \aleph_0$, then $Z$ is $\sigma$-closed discrete in $\bigcup \mc{V}$.}
\end{description}
\end{lem}

The conclusion of Lemma~\ref{todcount} had previously been shown under $\mr{MA}_{\omega_1}$ by Balogh~\cite{Balogh1983}. The weaker conclusion asserting that $Z$ is $\sigma$-discrete, if it's locally countable, was established by Todorcevic. A modification of his proof yields the stronger result~\cite{fischerTallNew}. It follows that:
\begin{cor}[]
\label{pfalc}
$\mr{PFA}(S)[S]$ implies that if $X$ is locally compact, includes no perfect pre-image of $\omega_1$, and $L(X)\leq\aleph_1$, and $Y\subseteq X$, $|Y| = \aleph_1$, is such that each point in $X$ has a neighbourhood meeting at most countably many points of $Y$, then $Y$ is $\sigma$-closed-discrete.
\end{cor}

We now need some results of Nyikos:
\begin{defn}[]
A space $X$ is of \textbf{Type I} if $X = \bigcup_{\alpha < \omega_1}X_{\alpha}$, where each $X_{\alpha}$ is open, $\alpha < \beta$ implies $\overline{X}_{\alpha}\subseteq X_{\beta}$, and each $\overline{X}_{\alpha}$ is Lindel\"{o}f. $\{ X_{\alpha} : \alpha < \omega_1\}$ is \textbf{canonical} if for limit $\alpha$, $X_{\alpha} = \bigcup_{\beta < \alpha}X_{\beta}$.
\end{defn}
\begin{lem}[~\cite{Nyikos2003}]
If $X$ is locally compact, $L(X) \leq\aleph_1$, and every Lindel\"{o}f subset of $X$ has Lindel\"{o}f closure, then $X$ is of Type I, with a canonical sequence.
\end{lem}
\begin{lem}[~\cite{Nyikos1983}]
If $X$ is of Type I, then $X$ is paracompact if and only if $\{\alpha : \overline{X}_{\alpha} - X_{\alpha}\neq 0\}$ is non-stationary.
\end{lem}
\begin{Prf}[\ref{nppthm}]
If $X$ is paracompact, this is straightforward. Suppose $X$ were not paracompact. $X$ is of Type I so we may pick a canonical sequence and  we may pick a stationary $S\subseteq\omega_1$ and $x_{\alpha}\in\overline{X}_{\alpha} - X_{\alpha}$, for each $\alpha\in S$. By Corollary~\ref{pfalc}, $\{ x_{\alpha} : \alpha\in S\}$ is $\sigma$-closed-discrete, so there is a stationary set of limit ordinals $S'\subseteq S$ such that $\{ x_{\alpha} : \alpha\in S'\}$ is closed discrete. Let $\{ U_{\alpha} : \alpha\in S'\}$ be a discrete collection of open sets expanding it. Pressing down yields an uncountable closed discrete subspace of some $X_{\alpha}$, contradiction.
\end{Prf}

Note that Lemma~\ref{lchncounter} follows from Theorem~\ref{mfpfa}, for consider the closure of a Lindel\"{o}f subspace $Y$ of a locally compact, hereditarily normal space which does not include a perfect pre-image of $\omega_1$. The following argument in Nyikos~\cite{Nyikos2003} will establish that $\overline{Y}$ is Lindel\"{o}f. First consider the special case when $Y$ is open. Let $B$ be a right-separated subspace of the boundary of $Y$. We claim $B$ is countable, whence the boundary is hereditarily Lindel\"{o}f, so $\overline{Y}$ is Lindel\"{o}f. Since the one-point compactification of $\overline{Y}$ is countably tight, by Lemma~\ref{todcount}, if $B$ is uncountable, it has a discrete subspace $D$ of size
$\aleph_1$. $D$ is closed discrete in $Z = \overline{Y} - (\overline{D} - D)$, so in $Z$ there is a discrete open expansion $\{ U_d : d\in D\}$ of $D$, because $\overline{Y}$ is hereditarily strongly $\aleph_1$-collectionwise Hausdorff by \textbf{CWH}. Since $Y\subseteq Z$, $\{ U_d \cap Y : d\in D\}$ is a discrete collection of non-empty subsets of $Y$, contradicting $Y$'s Lindel\"{o}fness.

Now consider an arbitrary Lindel\"{o}f $Y$. Since $X$ is locally compact, $Y$ can be covered by countably many open Lindel\"{o}f sets. The closure of their union is Lindel\"{o}f and includes $\overline{Y}$.

We next note that the requirement that Lindel\"{o}f subspaces have Lindel\"{o}f closures can be weakened. Recall the following result in~\cite{Gruenhage1996}:

\begin{lem}
Every locally compact, metalindel\"{o}f, $\aleph_1$-collectionwise Hausdorff, normal space is paracompact.
\end{lem}

Since in a normal $\aleph_1$-collectionwise Hausdorff space the closure of a Lindel\"{o}f subspace has countable extent, and metalindel\"{o}f spaces with countable extent are Lindel\"{o}f, we see that it suffices to have that closures of Lindel\"{o}f subspaces are metalindel\"{o}f.

\begin{cor}[]
There is a model of form $\mr{PFA}(S)[S]$ in which a locally compact normal space $X$ is paracompact if and only if the closure of every Lindel\"of subset of $X$ is metalindel\"of and $X$ does not include a perfect pre-image of $\omega_1$.
\end{cor}

A consequence of Corollary~\ref{pfalc} is that we can improve Theorem~\ref{nppthm} for spaces with Lindel\"of number $\leq \aleph_1$ to get:
\begin{thm}
\label{Lcor}
$\mr{PFA}(S)[S]$ implies that if $X$ is a locally compact normal space with $L(X)\leq\aleph_1$, and $X$ includes no perfect pre-image of $\omega_1$, then $X$ is paracompact.
\end{thm}
\begin{Proof}
As before, it suffices to consider the case of an open Lindel\"of $Y$. If the closure of $Y$ were not Lindel\"of, since it has Lindel\"of number $\leq \aleph_1$ there would be a locally countable subspace $Z$ of size $\aleph_1$ included in $\overline{Y} - Y$. That subspace would then be $\sigma$-closed-discrete by Corollary~\ref{pfalc}. As in the proof of Lemma~\ref{lchncounter} from Theorem~\ref{mfpfa}, we obtain a contradiction by getting an uncountable closed discrete subspace of $Y$. Since we have $\sigma$-\textit{closed}-discrete, we only need normality rather than hereditary normality.
\end{Proof}

In retrospect, Theorem~\ref{Lcor} is perhaps not so surprising: a phenomenon first evident in~\cite{Balogh1983} is that ``normal plus $L\leq\aleph_1$'' can often substitute for ``hereditarily normal'' in this area of investigation.

In fact, an even further weakening is possible:
\begin{defn}
Let $\mc{U}$ be an open cover of a space $X$ and let $x\in X$. $\mathbf{\op{\mathbf{Ord}}\boldsymbol{(}x,\mc{\mathbf{U}}\boldsymbol{)}} = |\{ U\in\mc{U} : x\in\mc{U}\}|$. $X$ is \textbf{submeta-$\boldsymbol{\aleph_1}$-Lindel\"{o}f} if every open cover has a refinement $\bigcup_{n < \omega}\mc{U}_n$ such that each $\mc{U}_n$ is an open cover, and for each $x\in X$, there is an $n$ such that $|\op{Ord}(x,\mc{U}_n)|\leq\aleph_1$.
\end{defn}
\begin{thm}
\label{smal}
There is a model of form $\mr{PFA}(S)[S]$ in which a locally compact normal space is paracompact and countably tight if and only if it is submeta-$\aleph_1$-Lindel\"{o}f and does not include a perfect pre-image of $\omega_1$.
\end{thm}
\begin{PROOF}
It suffices to prove that closures of Lindel\"{o}f subspaces have Lindel\"{o}f number $\leq\aleph_1$, for then we can apply Theorem~\ref{Lcor} to get that closures of Lindel\"{o}f subspaces are Lindel\"{o}f. Thus all we need is
\end{PROOF}
\begin{lem}[]
Every submeta-$\aleph_1$-Lindel\"{o}f space with extent $\leq\aleph_1$ has Lindel\"{o}f number $\leq\aleph_1$.
\end{lem}
\begin{Proof}
Following a similar proof in~\cite{Balogh1986}, suppose the space $X$ has an open cover $\mc{V}$ with no subcover of size $\leq\aleph_1$. Let $\bigcup_{n < \omega}\mc{U}_n$ be an open refinement of $V$ as in the definition of submeta-$\aleph_1$-lindel\"{o}fness. For each $y\in X$, pick $n(x)\in\omega$ such that $|\op{Ord}(x,\mc{U}_{n(x)})|\leq\aleph_1$. Let $X_n = \{ x\in X : n(x) = n\}$. Then for every $n < \omega$, there is a maximal $A_n\subseteq X_n$ such that no member of $\mc{U}_n$ contains two points of $A_n$. By maximality, $\mc{V}' = \bigcup_{n < \omega}\{\bigcup\{ U\in\mc{U}_n : x\in\mc{U}\} : x\in A_n\}$ covers $X_n$. Since $\mc{V}$ has no subcover of size $\leq\aleph_1$, $|A_n| > \aleph_1$ for some $n$. But $|A_n|$ is closed discrete.
\end{Proof}

An immediate corollary of Theorem~\ref{Lcor} is:
\begin{cor}[]
$\mr{PFA}(S)[S]$ implies every locally compact normal space of size $\leq\aleph_1$ with a $G_{\delta}$-diagonal is metrizable.
\end{cor}

The point is that spaces with a $G_{\delta}$-diagonal do not admit perfect pre-images of $\omega_1$, compact spaces with a $G_{\delta}$-diagonal are metrizable, and paracompact locally metrizable spaces are metrizable.

Weiss' space mentioned above constrains attempts at characterizing paracompactness. It is submetrizable, so has a $G_{\delta}$-diagonal. That latter property is hereditary; Lutzer~\cite{Lutzer1969} proved that linearly ordered spaces with a $G_{\delta}$-diagonal are metrizable, so:
\begin{prop}[]
If a space has a $G_{\delta}$-diagonal, it has no subspace homeomorphic to a stationary set.
\end{prop}

We are thus going to need stronger constraints on sets of size $\aleph_1$ than just excluding copies of stationary sets, if we wish to weaken the Lindel\"{o}f requirement of Theorem~\ref{mfpfa} to just something involving $\aleph_1$. Also note that Weiss' space prevents us from removing the cardinality restriction from Theorem~\ref{Lcor}. We will consider some such constraints in Section 4.

\section{Applications of P-ideal Dichotomy}
In order to prove Theorem \ref{thm2p6}, we introduce some known ideas about ideals.

\begin{defn}
A collection $\mc{I}$ of countable subsets of a set $X$ is a \textbf{$\mathbf{P}$-ideal} if each subset of a member of $\mc{I}$ is in $\mc{I}$, finite unions of members of $\mc{I}$ are in $\mc{I}$, and whenever $\{ I_n : n\in\omega\}\subseteq\mc{I}$, there is a $J\in\mc{I}$ such that $I_n - J$ is finite for all $n$.

$\mathbf{P}$ (short for \textbf{$\mathbf{P}$-ideal Dichotomy}): For every $P$-ideal $\mc{I}$ on a set  $X$, either
\begin{itemize}
\item[i)]there is an uncountable $A\subseteq X$ such that $[A]^{\leq\omega}\subseteq\mc{I}$
\item[or ii)]$X = \bigcup_{n < \omega}B_n$ such that for each $n$, $B_n\cap I$ is finite, for all $I\in\mc{I}$.
\end{itemize}
\end{defn}

\begin{defn}{\cite{EisworthNyikos}}
An ideal $\mc{I}$ of subsets of a set $X$ is \textbf{countable-covering} if for each countable $Q\subseteq X$, there is $\{ I_n^Q : n\in\omega\}\subseteq\mc{I}$, such that for each $I\in\mc{I}$ such that $I\subseteq Q$, $I\subseteq I_n^Q$ for some $n$.

$\mathbf{CC}$: For every countable-covering ideal $\mc{I}$ on a set $X$, either:
\begin{itemize}
\item[i)]there is an uncountable $A\subseteq X$ such that $A\cap I$ is finite, for all $I\in\mc{I}$,
\item[or ii)]$X = \bigcup_{n < \omega}B_n$ such that for each $n$, $[B_n]^{\leq\omega}\subseteq I$.
\end{itemize}
\end{defn}

Todorcevic's proof that $\mathop{PFA}(S)[S]$ implies $\mathbf{P}$ appears
in~\cite{LarsonErice}. In~\cite{EisworthNyikos}, Eisworth and Nyikos
proved that $\mathbf{P}$ implies $\mathbf{CC}$, and also proved the following remarkable result:

\begin{lem}
\label{CC}
$\mathbf{CC}$ implies that if $X$ is a locally compact space, then either
\begin{itemize}
\item[a)]$X$ is the union of countably many $\omega$-bounded subspaces,
\item[or b)]$X$ does not have countable extent,
\item[or c)]$X$ has a separable closed subspace which is not Lindel\"{o}f.
\end{itemize}
\end{lem}

Recall a space is \textbf{$\boldsymbol{\omega}$-bounded} if every countable subspace has compact closure. $\omega$-bounded spaces are obviously countably compact.

From~\cite{Gruenhage1977} we have:

\begin{lem}\label{obs}
An $\omega$-bounded space is either compact or includes a perfect pre-image of $\omega_1$.
\end{lem}

We can now prove Theorem \ref{thm2p6}.

The forward direction follows from \ref{mfpfa}.  To prove the other direction, it suffices to show that if $Y$ is a Lindel\"{o}f subspace of our space $X$, then $\overline{Y}$ is Lindel\"{o}f.  Applying \ref{CC}, we see that by \ref{obs}, $\overline{Y}$ will be $\sigma$-compact if we can exclude alternatives b) and c).  c) is excluded by hypothesis, so it suffices to show that $\overline{Y}$ has countable extent.  But that is easily established, since $\overline{Y}$ is locally compact normal and hence \textbf{$\aleph_1$-CWH}.  A closed discrete subspace of size $\aleph_1$ in $\overline{Y}$ could thus be fattened to a \textit{discrete} collection of open sets.  Their traces in $Y$ would contradict its Lindel\"{o}fness.
\hfill\qed.

\begin{cor}
There is a model of form $\mr{PFA}(S)[S]$ in which a locally compact space is metrizable if and only if it is normal, has a $G_{\delta}$-diagonal, and every separable closed subspace is Lindel\"{o}f.
\end{cor}
\begin{Proof}
Theorem~\ref{mfpfa} applies, since spaces with $G_{\delta}$-diagonals do not include perfect pre-images of $\omega_1$.
\end{Proof}

This characterization does not hold in $\mr{ZFC}$; the tree topology on a special Aronszajn tree is a locally compact Moore space, and hence has a $G_{\delta}$-diagonal. Under $\mr{MA}_{\omega_1}$, it is (hereditarily) normal. See e.g.\ the survey article~\cite{Tall1984}. Every separable subspace of an $\omega_1$-tree is bounded in height, and so is countable.

The condition in Theorem~\ref{thm2p6} that separable closed sets are Lindel\"{o}f, i.e. countable sets have Lindel\"of closures, can actually be weakened by a different argument, although perhaps the proof is more interesting than the result.

\begin{thm}\label{thm20}
There is a model of form $\mr{PFA}(S)[S]$ is which a locally compact normal space is paracompact and countably tight if and only if it includes no perfect pre-image of $\omega_1$ and the closure of each countable discrete subspace is Lindel\"of.
\end{thm}

We need several auxiliary results before proving this.

\begin{lem}[~\cite{ArhangBuzya1999}]\label{arhbuzlem}
If $X$ is Tychonoff, countably tight, $\aleph_1$-Lindel\"of, and countable discrete sets have Lindel\"of closures, then $X$ is Lindel\"of.
\end{lem}

Recall a space is defined to be \textbf{$\aleph_1$-Lindel\"of} if every open cover of size $\aleph_1$ has a countable subcover; equivalently, if every set of size $\aleph_1$ has a complete accumulation point.

\begin{thm}\label{thm25}
Assume $\mr{PFA}(S)[S]$.  Let $X$ be locally compact, normal, and not include a perfect pre-image of $\omega_1$.  Then either:
\begin{enumerate}
\item[a)] $X$ is $\sigma$-compact,
\item[or b)] $e(X) > \aleph_0$,
\item[or c)] $X$ has a countable discrete subspace $D$ such that $\overline D$ is not Lindel\"of.
\end{enumerate}
\end{thm}

\begin{Proof}
Assume b) and c) fail.  Since $X$ is locally compact and countably tight, by Lemma~\ref{arhbuzlem}, it suffices to prove $X$ is $\aleph_1$-Lindel\"of.

If not, there is a $Y \subseteq X$ of size $\aleph_1$ with no complete accumulation point.  Thus $Y$ is locally countable and hence $\sigma$-discrete.  Hence there is an uncountable discrete $Z \subseteq Y$ with no complete accumulation point.  Let $Z = \{z_\alpha:\alpha<\omega_1\}$.  Then $\overline{Z} = \bigcup_{\beta < \omega_1}\overline{\{z_\alpha:\alpha<\beta\}}$.  By hypothesis, it follows that $L(\overline{Z})\leq\aleph_1$.  But then $\overline{Z}$ is paracompact by Theorem \ref{mfpfa}.  But $e(X)\leq\aleph_0$, so $\overline{Z}$ is Lindel\"of, so $Z$ does have a complete accumulation point, giving a contradiction.
\end{Proof}

\begin{trivlist} \item[] {\bf Proof of Theorem~\ref{thm20}.}
It suffices to show closures of countable subspaces of our spaces are Lindel\"of.  By normality and $\aleph_1$-collectionwise Hausdorfness they have countable extent, and we are assuming countable discrete subspaces have Lindel\"of closures, so we are done by Theorem~\ref{thm25}
\hspace*{0pt}\hfill $\Box$\end{trivlist}

One is tempted to substitute the condition that \textit{discrete subspaces have Lindel\"of closures} for b) and c) of Theorem~\ref{thm25}, but unfortunately it is already known that countably tight spaces with that property are Lindel\"of~\cite{ArhangBuzya1999}.

\section{Applications of Axiom R and Hereditary Paracompactness}
We shall consider hereditary paracompactness and obtain some interesting results. We will need  the following result of Balogh~\cite{Balogh2003}:
\begin{lem}[]
\label{ctctr}
If $X$ is countably tight, has a dense subspace of size $\leq\aleph_1$, and every subspace of size $\leq\aleph_1$ is metalindel\"{o}f, then $X$ is hereditarily metalindel\"{o}f.
\end{lem}

Balogh assumes in \cite{Balogh2003} that all spaces considered are regular, but does not use regularity in the proof of Lemma \ref{ctctr}.  He also does not actually require \textit{all} subspaces of size $\leq\aleph_1$ to be meta\-lindel\"{o}f in order to obtain the conclusion of Lemma~\ref{ctctr}. We refer the reader to~\cite{Balogh2003} for the details. Note that it follows that in a countably tight space in which every subspace of size $\leq\aleph_1$ is metalindel\"{o}f, separable sets are (hereditarily) Lindel\"{o}f. Also note that Weiss' space must have a subspace of size $\aleph_1$ which is not metalindel\"{o}f.

\begin{thm}
\label{axr}
Axiom $R$ implies a locally separable, regular, countably tight space is hereditarily paracompact if and only if every subspace of size $\leq\aleph_1$ is meta\-lindel\"{o}f.
\end{thm}
\begin{Proof}
One direction is trivial. To go the other way, we shall first obtain paracompactness via Lemma~\ref{rimp}. Here we do need regularity.  I thank Sakae Fuchino for pointing this out.  Let $V$ be an open subspace with $L(V)\leq\aleph_1$. Covering $V$ by $\leq\aleph_1$ separable open sets, we see that $d(V)\leq\aleph_1$. Then by Lemma~\ref{ctctr}, $\overline{V}$ is hereditarily paracompact. To get the whole space hereditarily paracompact, note it is a sum of separable, hence hereditarily Lindel\"{o}f, clopen sets.
\end{Proof}

\begin{cor}[]
Axiom $R$ implies that a locally hereditarily separable, regular space is hereditarily paracompact if and only if each subspace of size $\leq\aleph_1$ is metalindel\"{o}f.
\end{cor}

\begin{Proof}
Local hereditary separability implies countable tightness.
\end{Proof}

\begin{cor}[]
Axiom $R$ implies a locally second countable, regular space is metrizable if and only if every subspace of size $\leq\aleph_1$ is metalindel\"{o}f.
\end{cor}
\begin{Proof}
This is clear, since such a space is locally hereditarily separable, while paracompact, locally metrizable spaces are metrizable.
\end{Proof}

\begin{cor}[~\cite{Balogh2002}]
Axiom $R$ implies every locally compact space in which every subspace of size $\aleph_1$ has a point-countable base is metrizable.
\end{cor}
\begin{Proof}
Dow~\cite{Dow1988} showed that compact spaces in which every subspace of size $\aleph_1$ has a point-countable base are metrizable.
\end{Proof}

\begin{cor}[]
Axiom $R$ implies a locally compact space is metrizable if and only if it has a $G_{\delta}$-diagonal and every subspace of size $\leq\aleph_1$ is metalindel\"{o}f.
\end{cor}
\begin{Proof}
Compact spaces with $G_{\delta}$-diagonals are metrizable.
\end{Proof}
\noindent\textbf{Note: }Results similar to ours concerning Axiom R were obtained by S. Fuchino and his collaborators independently \cite{Fuchino}, \cite{Fuchinoa}.

With the added power of PFA$(S)[S]$, we can utilize Lemma \ref{ctctr} without assuming local separability.  First, we observe:
\begin{thm}
\label{xct}
If $X$ is countably tight and every subspace of size $\leq\aleph_1$ is metalindel\"{o}f, then $X$ does not include a perfect pre-image of $\omega_1$.
\end{thm}

\begin{prf}[\ref{xct}]
This follows immediately from Lemma~\ref{ctctr} and:

\begin{lem}\label{lem4p10}
Every perfect pre-image of $\omega_1$ includes one of density $\leq\aleph_1$.
\end{lem}

\begin{Proof}
  Let $\pi : X\rightarrow\omega_1$ be perfect and onto. Let $C =
  \{\alpha : \pi^{-1}(\alpha + 1) - \overline{\pi^{-1}(\alpha)}\neq
  0\}$. Then $C$ is unbounded, for suppose not. Then there is an
  $\alpha_0$ such that $Y = \pi^{-1}(\alpha + 1)$. But then $Y$ is
  compact, contradiction. Pick for each $\alpha\in C$, a
  $d_{\alpha}\in\pi^{-1}(\alpha + 1) -
  \overline{\pi^{-1}(\alpha)}$. Let $Q = \overline{\{ d_{\alpha} :
    \alpha\in C\}}$. Then $\pi\mid Q$ is perfect, so $\pi(Q)$ is
  closed unbounded, so is homeomorphic to $\omega_1$.
\end{Proof}
\vspace*{-10pt}
\end{prf}
From Theorem \ref{xct} we then obtain:

\begin{thm}
In the model of form $\mr{PFA}(S)[S]$ of \cite{Larson} a locally compact, normal, countably tight space is paracompact if every subspace of size $\aleph_1$ is metalindel\"{o}f.
\end{thm}

\begin{Proof}
Separable subspaces are Lindel\"{o}f; by \ref{lem4p10} and \ref{ctctr}, there are no perfect pre-images of $\omega_1$ in such spaces.
\end{Proof}

\begin{cor}
\label{mpfass}
In the $\mr{PFA}(S)[S]$ model of \cite{Larson}, a locally compact, countably tight space is hereditarily paracompact if and only if it is hereditarily normal and every subspace of size $\leq\aleph_1$ is metalindel\"{o}f.
\end{cor}
\begin{Proof}
This follows from Theorem \ref{xct} and the observation that countable tightness is inherited by open subspaces.  
\end{Proof}

\begin{cor}\label{lclss}
There is a model of form $\mathrm{PFA}(S)[S]$ in which a locally compact, locally separable space is hereditarily paracompact if and only if it is hereditarily normal and every subspace of size $\leq\aleph_1$ is metalindel\"{o}f.
\end{cor}
This will follow from Theorem~\ref{mpfass} and \textbf{FC $\aleph_1$-CWH}, since the latter implies the hypothesis of the following:

\begin{lem}[~\cite{Nyikos1992}]\label{snfcs}
If separable, normal, first countable spaces do not have uncountable closed discrete subspaces, then compact, separable, hereditarily normal spaces are countably tight.
\end{lem}

We can avoid introducing the hitherto unused axiom \textbf{FC $\aleph_1$-CWH} by quoting:

\begin{lem}[~\cite{Tall1977}]\label{lem4p13}
If there is a separable, normal, first countable space with an uncountable closed discrete subspace, there is a locally compact one.
\end{lem}

\begin{proofcor}[\ref{lclss}]
It suffices to show such a space is countably tight. Given $x\in\overline{Y}$, there is a separable open neighbourhood $U$ of $x$ with $\overline{U}$ compact. Then $x\in\overline{(U\cap Y)}$. $\overline{U}$ is countably tight by Lemma~\ref{snfcs}. Thus there is a countable $D\subseteq U\cap Y$ such that $x\in\overline{D}\cap\overline{U}$. But then $x\in\overline{D}$ as required.
\end{proofcor}

There is another way of proving Corollary~\ref{lclss}, which actually gives a slightly stronger result: locally satisfying the countable chain condition instead of locally separable. This follows from~\cite{Todorcevic}, in which Todorcevic showed that \textit{$\mathit{PFA}$ implies compact, hereditarily normal spaces satisfying the countable chain condition are hereditarily Lindel\"{o}f} (and hence first countable).  Since \cite{Todorcevic} is still unavailable, we shall prove this in Section 6.

\section{(Local)\-Connectedness and ZFC Reflections}
One can sometimes replace our use of Axiom $\mr{R}$ by the assumption of (local) connectedness, thanks to the following observation:

\begin{lem}[{~\cite[5.9]{EisworthNyikos}}]
\label{xlcct}
If $X$ is locally compact, locally connected, and countably tight, then $X$ is a topological sum of Type I spaces if and only if every Lindel\"{o}f subspace of $X$ has Lindel\"{o}f closure.  Similarly, if $X$ is locally compact, connected, countably tight, and Lindel\"{o}f subspaces have Lindel\"{o}f closures, then it is Type I.
\end{lem}

Thus we have:

%2.12'%
\setcounter{section}{2}
\setcounter{theorem}{11}
\begin{thmp}\label{2p12p}
$\mr{PFA}(S)[S]$ implies a locally compact, locally connected, normal space $X$ is paracompact if and only if separable closed subspaces are Lindel\"{o}f, and $X$ does not include a perfect pre-image of $\omega_1$.
\end{thmp}

%4.2'%
\setcounter{section}{4}
\setcounter{theorem}{1}
\begin{thmp}
\label{lclclsct}
A locally compact, (locally) connected, locally separable, countably tight, regular space is hereditarily paracompact if and only if every subspace of size $\leq\aleph_1$ is metalindel\"{o}f.
\end{thmp}

\begin{Proofs}
  Theorem~\ref{lclclsct}$'$ is the only one which requires a bit of
  thought. Any Lindel\"of subspace is included in a separable open set
  $S$. $\overline{S}$ is Lindel\"of and therefore so is
  $\overline{L}$. Thus the space is a sum of Type I spaces, each of
  density $\leq\aleph_1$, and by Lemma~\ref{ctctr}, each of these is
  hereditarily metalindel\"{o}f. By local separability, the space is
  then hereditarily paracompact.
\end{Proofs}

%4.3'%
\begin{corp}\label{lclclhss}
A locally compact, (locally) connected, locally hereditarily separable, regular space is hereditarily paracompact if and only if each subspace of size $\leq\aleph_1$ is metalindel\"{o}f.
\end{corp}

Particularly pleasant is:

\setcounter{section}{5}
\setcounter{theorem}{1}
\begin{cor}
A manifold is metrizable if and only if every subspace of size $\aleph_1$ is metalindel\"{o}f.
\end{cor}

\setcounter{section}{4}
\setcounter{theorem}{4}
%4.5'%
\begin{corp}
A locally compact, (locally) connected space in which every subspace of size $\aleph_1$ has a point-countable base is metrizable.
\end{corp}

%4.6'%
\begin{corp}
A locally compact, (locally) connected space is metrizable if and only if it has a $G_{\delta}$-diagonal and every subspace of size $\aleph_1$ is metalindel\"{o}f.
\end{corp}

We also have:
\setcounter{theorem}{6}
\begin{thmp}
$\mr{PFA}(S)[S]$ implies a locally compact, normal, countably tight, connected or locally connected space is paracompact if every subspace of size $\aleph_1$ is metalindel\"{o}f.
\end{thmp}
\setcounter{theorem}{8}
%4.9'%
\begin{thmp}
$\mr{PFA}(S)[S]$ implies a locally compact, locally connected, countably tight space is hereditarily paracompact if and only if it is hereditarily normal and every subspace of size $\leq\aleph_1$ is metalindel\"{o}f.
\end{thmp}

\setcounter{theorem}{11}
%4.12'%
\begin{corp}\label{4p12p}
$\mr{PFA}(S)[S]$ implies a locally compact, locally connected, locally separable space is hereditarily paracompact if and only if it is hereditarily normal and every subspace of size $\leq\aleph_1$ is metalindel\"{o}f.
\end{corp}

Balogh~\cite{Balogh2002} proved:

\setcounter{section}{5}
\setcounter{theorem}{2}
\begin{lem}
\label{xllct}
Let $X$ be a locally Lindel\"{o}f, regular, countably tight space with $L(X)\leq\aleph_1$. Suppose that every subspace of size $\leq\aleph_1$ of $X$ is paracompact, and $X$ is either normal or locally has countable spread. Then $X$ is paracompact.
\end{lem}

We then have the following variation of Corollary~\ref{lclclhss}$'$:

\begin{thm}
\label{xlclcs}
Let $X$ be a locally compact, (locally) connected space in which every subspace of size $\leq\aleph_1$ is metalindel\"{o}f, and which locally has countable spread. Then $X$ is hereditarily paracompact.
\end{thm}

\begin{PROOF}
  It suffices to show $X$ is paracompact, since all the properties in
  question are open-hereditary. By Lemmas~\ref{xlcct} and~\ref{xllct},
  it suffices to prove that $X$ is countably tight and closures of Lindel\"{o}f subspaces are
  Lindel\"{o}f. Lindel\"{o}f subspaces are included in the union of
  countably many subspaces with countable spread and hence have
  countable spread. If a Lindel\"{o}f $Y\subseteq X$ did not have
  (hereditarily) Lindel\"{o}f closure, there would be a
  right-separated subset $Z$ of $\overline{Y}$, with $|Z| =
  \aleph_1$. But $Z$ would then be metalindel\"{o}f and locally
  countable, hence paracompact and $\sigma$-discrete. Note that $X$
  --- and hence $\overline{Z}$ --- is countably tight, since compact
  spaces with countable spread are countably
  tight~\cite{Arhangelskii1971}. By Lemma~\ref{ctctr}, $\overline{Z}$
  is then hereditarily metalindel\"{o}f. $\overline{Z}$ is locally
  separable, since if $U$ is an open subspace of $X$ with countable
  spread, $Z \, \cap \, U$ is dense in $\overline{Z} \, \cap \, U$,
  but is countable, since $Z$ is $\sigma$-discrete. Similarly the
  closure of $Z$ in $Y$ is locally separable. But the closure of $Z$
  in $Y$ is Lindel\"{o}f, so it's separable. But then $\overline{Z}$
  is separable. But then $\overline{Z}$ is hereditarily Lindel\"{o}f,
  contradiction.
\end{PROOF}

The advantage of eliminating explicit and implicit uses of Axiom $R$ as we did in \ref{2p12p}$'$  and \ref{4p12p}$'$ is that it makes it likely that such results can then be obtained without the necessity of assuming large cardinals, by using $\aleph_2$-p.i.c.\ forcing as in e.g.~\cite{Todorcevic1985}.

\section{PFA$(S)[S]$ and Locally Compact Dowker Spaces}
The question of whether there exist \textit{small Dowker spaces}, i.e.\ normal spaces with product with the unit interval not normal, which have familiar cardinal invariants of size $\leq 2^{\aleph_0}$, continues to attract attention from set-theoretic topologists. See for example the surveys~\cite{Balogh2004,Rudin1984,Szeptycki,Szeptycki1993}. Although there are many consistent examples, there have been very few results asserting the consistency of the non-existence of such examples. We shall partially remedy that situation here.  In this section, we observe that $\mr{PFA}(S)[S]$ excludes some possible candidates for \textit{small Dowker spaces}. Most of our results follow easily from what we have already proved. Recall:

\begin{lem}[~\cite{Dowker1951}]
\label{nsx}
For a normal space $X$, the following are equivalent:
\begin{itemize}
\item[a)]$X$ is countably paracompact,
\item[b)]$X\times [0,1]$ is normal,
\item[c)]$X\times (\omega+1)$ is normal.
\end{itemize}
\end{lem}

``Small'' is not very well-defined; in the recent survey~\cite{Szeptycki}, Szeptycki concentrates on the properties \textit{cardinality $\aleph_1$}, \textit{first countability}, \textit{separability}, \textit{local compactness}, \textit{local countability} (i.e.\ each point has a countable neighbourhood) and \textit{submetrizability} (i.e.\ the space has a weaker metrizable topology). We shall deal with several of these, weakening --- in terms of non-existence --- \textit{cardinality $\leq\aleph_1$} to \textit{Lindel\"{o}f number $\leq\aleph_1$} and \textit{submetrizable} to \textit{not including a perfect pre-image of $\omega_1$}. Note that submetrizable spaces have $G_{\delta}$-diagonals and hence cannot include perfect pre-images of $\omega_1$, since countably compact spaces with $G_{\delta}$-diagonals  are metrizable~\cite{Chaber1976}.

\begin{thmm}\mbox{}
\begin{enumerate}
\item[1)]$\mr{PFA}(S)[S]$ implies there is no locally compact, hereditarily normal Dowker space which in addition:
\begin{itemize}
\item[a)]satisfies the countable chain condition,
\item[or b)]includes no perfect pre-image of $\omega_1$ and is either
  connected or locally connected.
\item[or c)]has countable extent.
\end{itemize}
\item[2)]$\mr{PFA}(S)[S]$ implies there is no locally compact Dowker space which includes no perfect pre-image of $\omega_1$ and has Lindel\"{o}f number $\leq\aleph_1$.
\item[3)]In the PFA$(S)[S]$ model of \cite{Larson}:
\begin{enumerate}
\item there is no locally compact, hereditarily normal Dowker space including a perfect pre-image of $\omega_1$.
\item there is no locally compact Dowker space in which separable closed subspaces are Lindel\"{o}f and which includes no perfect pre-image of $\omega_1$.
\item there is no locally compact, countably tight Dowker space in which every subspace of size $\aleph_1$ is \mbox{meta\-lindel\"{o}f}.
\item there is no locally compact, countably tight, Dowker $D$-space.
%\item[7)]$\mr{PFA}(S)[S]$ implies there is no locally compact, discretely Lindel\"of Dowker space.
%\marginpar{I can move ``meta-'' down if you'd prefer}
\end{enumerate}
\end{enumerate}
\end{thmm}

We shall start with:
\begin{thm}
\label{shn_thm}
Assume $\mr{PFA}(S)[S]$. Let $X$ be a locally compact, hereditarily normal space satisfying the countable chain condition. Then $X$ is hereditarily Lindel\"{o}f, and hence countably paracompact.
\end{thm}

\begin{Proof}
This follows from \cite{Todorcevic}, where Todorcevic proves:
\begin{lem}
\label{pfalem}
$\mr{PFA}(S)[S]$ implies compact hereditarily normal spaces satisfying
the countable chain condition are hereditarily Lindel\"{o}f.
\end{lem}
\begin{Proof}
Since open subspaces are locally compact normal, the space is hereditarily $\aleph_1$-collectionwise-Hausdorff and hence has countable spread.  If the space were not hereditarily Lindel\"{o}f, it would have an uncountable right-separated subspace, and hence, by $\boldsigma$, an uncountable discrete subspace, contradiction.
\end{Proof}

Todorcevic's proof was more difficult, since \textbf{$\aleph_1$-CWH} was not available to him.

The one-point compactification of a locally compact, hereditarily
normal space $X$ is hereditarily normal, and satisfies the
countable chain condition if and only if $X$ does. The result follows, so we have established 1a) of the Main Theorem.
\end{Proof}

1b) of the Main Theorem follows from \ref{Lcor} plus \ref{xlcct}.

To prove 1c), we call on Theorem \ref{thm2p6}.  Since separable closed subspaces are Lindel\"{o}f, the space is the union of countably many $\omega$-bounded -- hence countably compact -- subspaces.  In a normal space, the closure of a countably compact subspace is countably compact, and it is not hard to show that the union of countably many countably compact closed subspaces of a normal space is countably paracompact.

Restating 2) of the Main Theorem, we next have:

\begin{thm}\label{nolcds}
  $\mr{PFA}(S)[S]$ implies every locally compact Dowker space of
  Lindel\"{o}f number $\le \aleph_1$ includes a perfect
  pre-image of $\omega_1$.
\end{thm}

\begin{cor}\label{nolcdscor}
PFA$(S)[S]$ implies there are no locally compact submetrizable Dowker spaces of size $\aleph_1$.
\end{cor}

\begin{Proof}
\ref{nolcds} follows immediately from \ref{Lcor}.
\end{Proof}
The conclusion of Theorem~\ref{nolcds} is an improvement of a result in~\cite{Balogh1983}; Balogh proved from $\mr{MA}_{\omega_1}$ that locally compact spaces of size $\aleph_1$ which don't include a perfect pre-image of $\omega_1$ are $\sigma$-closed-discrete, hence, if normal, are countably paracompact.

To show that Theorem~\ref{shn_thm}  is not vacuous, we note that Nyikos~\cite{Nyikos1999} constructed, assuming $\largediamond$, a hereditarily separable,
locally compact, first countable, hereditarily normal Dowker space.

In~\cite{Juhasz1976}, the authors remark that they can construct under
$\largediamond$, using their technique of refining the topology on a
subspace of the real line, a locally compact Dowker space. By
$\mr{CH}$, such a space has cardinality $\aleph_1$. Since it refines
the topology on a subspace of $\mathbb{R}$, it is submetrizable. Thus
the conclusion of Corollary~\ref{nolcdscor} is independent. We do not have
consistent counterexamples for clauses 1b), 3b), c), d) of the Main
Theorem.

Clause 3a) of the Main Theorem follows immediately from \ref{lchncounter}; 3b) follows from \ref{thm2p6}.  To prove 3c), first observe that by \ref{xct}, $X$ does not include a perfect pre-image of $\omega_1$.  Next, if $Y$ is a separable closed subspace of $X$, by \ref{ctctr} $Y$ is Lindel\"{o}f.\hfill\qed

``$D$-spaces'' are popular these days.  See e.g.~\cite{Eisworth2007}, \cite{Gruenhage2009}.

\begin{defn}
  X is a \textbf{$D$-space} if for every neighborhood assignment
  $\{V_x\}_{x \in X}$, there is a closed discrete $Y \subseteq X$ such
  that $\bigcup \{V_x: x \in Y\}$ is a cover.
\end{defn}

\begin{thm}\label{thm19a}
There is a model of form $\mr{PFA}(S)[S]$ in which a locally compact normal countably tight space is paracompact if and only if it is a $D$-space.
\end{thm}

Clause 3d) of the Main Theorem follows. Theorem~\ref{thm19a} is
analogous to the fact that \textit{linearly ordered spaces are
  paracompact if and only if they are $D$-spaces} (see
e.g. \cite{Gruenhage2009}).

\begin{trivlist} \item[] {\bf Proof of Theorem~\ref{thm19a}.}
Assume the space is $D$.  It is well-known and easy to see that countably compact $D$-spaces are compact. It is also easy to see that closed subspaces of $D$-spaces are $D$.

It follows from Lemma \ref{lem2p4} that a countably tight $D$-space cannot include a perfect
pre-image of $\omega_1$. By \textbf{$\aleph_1$-CWH}, the closure of a
countable subspace of our space is collectionwise Hausdorff, and
hence has countable extent. But again, it is well-known that
$D$-spaces with countable extent are Lindel\"of. By Theorem \ref{thm2p6}, our space is then paracompact.

  For the other direction, a paracompact, locally compact space is a
  discrete sum of $\sigma$-compact spaces. It is well-known that
  $\sigma$-compact spaces are $D$-spaces, and it is easy to verify
  that discrete sums of $D$-spaces are $D$-spaces.
  \hspace*{0pt}\hfill $\Box$\end{trivlist}

One way of strengthening normality without necessarily implying countable paracompactness is to assume hereditary normality. Another is to assume powers of the space are normal. And then one could assume both. Let's see what happens. We have already looked at hereditary normality; but let us also recall from~\cite{Larson} that:

\begin{prop}
There is a model of form $\mr{PFA}(S)[S]$ in which every locally compact space with hereditarily normal square is metrizable.
\end{prop}

Even in $\mr{ZFC}$, a hereditarily normal square has consequences. The following results are due to P.\ Szeptycki~\cite{Szeptyckia}:

\begin{prop}
\label{szthm}
If $X^2$ is normal and $X$ includes a countable non-discrete subspace, then $X$ is countably paracompact.
\end{prop}

\begin{cor}
If $X$ is separable, first countable, or locally compact, and $X^2$ is normal, then $X$ is countably paracompact.
\end{cor}

On the other hand, following a suggestion of W.\ Weiss, Szeptycki \cite{Szeptyckia} noticed that Rudin's $\mr{ZFC}$ Dowker space~\cite{Rudin1971/72} has all finite products normal.

Although our consistency results concerning small Dowker spaces improve previous ones, they have two unsatisfactory aspects. First of all, all but 1c) prove paracompactness, rather than countable paracompactness, so there ought to be sharper results.

It is likely that in our results involving hereditary normality, ``perfect pre-image of $\omega_1$'' can be weakened to ``copy of $\omega_1$.'' This would follow from the following conjecture and unpublished theorem of the author.

\begin{con}
$\mr{PFA}(S)[S]$ implies every first countable perfect pre-image of $\omega_1$ includes a copy of $\omega_1$.
\end{con}

\begin{thm}
\label{ite}
$\mr{PFA}(S)[S]$ implies that every hereditarily normal perfect pre-image of
$\omega_1$ includes a first countable perfect pre-image of $\omega_1$.
\end{thm}

Another unsatisfactory aspect of our consistency results is that a
supercompact cardinal is required to construct models of form
$\mr{PFA}(S)[S]$. This is surely overkill, when we are really
concerned with $\aleph_1$. We suspect that large cardinals are not
needed except possibly for those relying on Axiom R.  The other
clauses probably can be obtained without any large cardinals, by $\aleph_2$-p.i.c.\ forcing as in e.g.~\cite{Todorcevic1985}.

\section{Hereditarily Normal Compact Spaces}
Under PFA$(S)[S]$, hereditarily normal compact spaces -- ``$T_5$ compacta'' for short -- have strong properties.  We have already seen (\ref{shn_thm}) that separable ones are hereditarily Lindel\"{o}f.  It follows that they are first countable.  Hence:
\begin{thm}\label{thm7p1new}
Countably compact, locally compact $T_5$ spaces are sequentially compact.
\end{thm}
\begin{cor}\label{thm7p1}
PFA$(S)[S]$ implies $T_5$ compacta are sequentially compact.
\end{cor}
\begin{Prf}[\ref{thm7p1new}]
Let $X$ be a countably compact, locally compact $T_5$ space.  The one-point compactification of the closure of the range of a sequence is a separable $T_5$ compactum, so is first countable.  The closure of the range is then itself first countable, so there is a subsequence converging to a limit point of the range.
\end{Prf}

Juh\'{a}sz, Nyikos, Szentmikl\'{o}ssy \cite{JNS} proved:

\begin{lem}\label{lem7p2}
$T_5$ compacta which are homogeneous and hereditarily strongly $\aleph_1$-collectionwise-Hausdorff are countably tight.
\end{lem}

It follows by \textbf{$\aleph_1$-CWH} that PFA$(S)[S]$ implies homogeneous $T_5$ compacta are countably tight. But we can do better:

\begin{thm}\label{thm7p3}
PFA$(S)[S]$ implies homogeneous $T_5$ compacta are first countable.
\end{thm}
\begin{Proof}
The authors of \cite{JNS} show that homogeneous $T_5$ compacta are first countable, provided their open Lindel\"{o}f subspaces have hereditarily Lindel\"{o}f boundaries.  We proved this following the proof of \ref{nppthm} above.
\end{Proof}

The conclusion of \ref{thm7p3} was earlier proved consistent by de la Vega, using a different model \cite{V}.

The conclusion of \ref{thm7p1} is not true in ZFC: Fedorchuk's $S$-space from $\diamondsuit$ \cite{Fe} is a $T_5$ compactum which is countably tight -- because it is hereditarily separable -- but has no non-trivial convergent sequences.

\begin{rem}
  The proofs in this paper were produced around 2006-2007, assuming \textbf{$\aleph_1$-CWH} for the PFA$(S)[S]$ ones. However, a correct proof
  of that was only obtained in 2010. At the 2006 Prague
  Topological Symposium, Todorcevic announced the $\sigma$-discrete
  version of $\boldsigma$ followed from PFA$(S)[S]$.  Larson \cite{LarsonErice} wrote some
  notes on Todorcevic's lectures at the conference on \emph{Advances in Set-theoretic
    Topology, in Honor of T. Nogura} in Erice, Italy in 2008
  \cite{Todorcevic2008}.  Using these and ideas of Todorcevic, A. Fischer and the author derived a proof of $\boldsigma$ from PFA$(S)[S]$ \cite{fischerTallNew}.
\end{rem}

\begin{ackn}
I am grateful to members of the Toronto Set Theory Seminar and to Gary Gruenhage and Sakae Fuchino for discussions concerning this work.  I also thank the referee of \cite{TallDowker} for correcting the statements of the $D$-space results, and for suggesting I merge \cite{TallDowker} and \cite{T1}.
\end{ackn}

\nocite{*}
\bibliographystyle{acm}
\bibliography{CP2}

{\rm Franklin D. Tall, Department of Mathematics, University of
Toronto, Toronto, Ontario M5S 2E4, CANADA}

{\it e-mail address:} {\rm tall@@math.utoronto.ca}

\end{document}